\documentclass[12pt]{article}
\usepackage{amssymb,amsfonts,amsthm,amsmath}
\usepackage{slashed} 
\usepackage{yfonts}
\usepackage{mathrsfs,pifont}

\usepackage[all]{xy}

\usepackage{amssymb,amsfonts,amsthm,amsmath}
\usepackage[all]{xy}
\usepackage{slashed} 
\usepackage{yfonts}
\usepackage{mathrsfs,pifont}

\newcommand{\C}{\mathbb{C}}
\newcommand{\Z}{\mathbb{Z}}

\newcommand{\Q}{\mathbb{Q}}

\newcommand{\qMZV}{\mathsf{qMZV}}
\newcommand{\grqMZV}{\mathsf{grqMZV}}
\newcommand{\MZV}{\mathsf{MZV}}
\newcommand{\QM}{\mathsf{QM}}
\newcommand{\oZ}{\mathsf{oZ}}
\newcommand{\oqZ}{\mathsf{oqZ}}

\newcommand{\bs}{\mathsf{s}}
\newcommand{\pt}{\textup{pt}}

\newcommand{\fS}{\mathfrak{S}}

\newcommand{\cK}{\mathscr{K}}

\newcommand{\cO}{\mathscr{O}}
\newcommand{\cL}{\mathscr{L}}
\newcommand{\cT}{\mathscr{T}}
\newcommand{\cI}{\mathscr{I}}
\newcommand{\cM}{\mathscr{M}}

\newcommand{\td}{\textup{\sc 3D}}
\newcommand{\vir}{\textup{vir}}
\newcommand{\bd}{\boldsymbol{\delta}}
\newcommand{\bdt}{\bd_\td}

\newcommand{\Ct}{{\C^{\times}}}

\newcommand{\lla}{\left\langle}
\newcommand{\rra}{\right\rangle^{\! \circ}}

\newcommand{\bx}{\textup{\ding{114}}}

\newcommand{\Mgn}{\overline{\mathcal{M}}_{g,n}}

\DeclareMathOperator{\weight}{weight}
\DeclareMathOperator{\depth}{depth}
\DeclareMathOperator{\Hilb}{Hilb}
\DeclareMathOperator{\Ext}{Ext}
\DeclareMathOperator{\Pic}{Pic}
\DeclareMathOperator{\Euler}{Euler}
\DeclareMathOperator{\ch}{ch}

\newtheorem{Conjecture}{Conjecture}

\begin{document}

\title{Hilbert schemes and multiple $q$-zeta values} 
\author{Andrei Okounkov} 
\date{} 
\maketitle

\begin{abstract}
We present several conjectures on multiple $q$-zeta values 
and on the role they play in certain problems of enumerative geometry. 
\end{abstract}

\section{Multiple $q$-zeta values}\label{s1} 

\subsection{}

Multiple zeta values, that is, series of the form
\begin{equation}
\zeta(\bs) = \zeta(s_1,\dots,s_k) = \sum_{n_1 > n_2 > \dots > n_k} 
\frac{1}{n_1^{s_1} \dots n_k^{s_k}}\,,
\label{MZV}
\end{equation}
where $s_i$ are nonnegative integers and $s_1>1$, 
drew the attention of some of the greatest mind in mathematics, starting 
with Euler, see e.g.\ \cite{Zag}. 
Many remarkable properties, striking 
applications, and unexpected connections were found for them, 
some conjectural, some
proven. For an introduction to their motivic and 
representation-theoretic aspects, the reader may turn, 
for example, to chapter 
25 in \cite{Andre} and the Appendix in \cite{ES}, respectively. 

\subsection{} 

In this paper, we work with $q$-deformations of the series
\eqref{MZV}. There already exists a large body of work on various 
flavors of such $q$-deformations, see e.g.\ \cite{Brad,OT,Zud}. 
The one we use here 
has the form
\begin{equation}
  \label{qMZV}
Z(\bs) = \sum_{n_1 > n_2 > \dots > n_k} 
\prod (n_i)^{-s_i}\,, \quad  (n)^{-s} = \frac{p_s(q^n)}{(1-q^n)^s}\,, 
\end{equation}
where $p_s$ is a nonzero polynomial of degree $s$ without constant
term. 

The linear span of the series \eqref{qMZV} is independent of the 
choice of the numerators $p_s$. However, to discuss a conjectural 
$\Z/2$-grading, it is convenient to choose the
numerators
palindromic, that is, satisfying $t^s p_s(1/t) = p_s(t)$, which is
only possible if $s>1$. In what follows, we assume that $s_i\ge 2$ 
and define 
$$
p_s(t) = 
\begin{cases}
t^{s/2}\,,  &s=2,4,6,\dots\,,
\\ 
t^{(s-1)/2}(1+t)\,,  &s=3,5,7,\dots\,. 
\end{cases}
$$

\subsection{}

We define $\qMZV$ as the $\Q$-subalgebra of $\Q[[q]]$ spanned  
by the series $Z(\bs)$ with $s_i \ge 2$. This algebra is 
filtered by weight, where 
$$
\weight \, Z(\bs)= \sum s_i \,.
$$
Note that 
$$
\QM= \Q[Z(2),Z(4),Z(6)] \subset \qMZV
$$
is the classical ring of quasimodular forms \cite{KZ}, the 
analog of $\Q[\pi^2] \subset \MZV$.  It contains all $Z(2k)$ 
and with proper normalization (which includes 
correct constant terms), is $\Z$-graded by weight.  
I don't know whether there is a way to 
extend this grading to all of $\qMZV$.

Clearly, 
$$
(1-q)^{\weight} \, Z(\bs) \to 2^{\# \textup{odd}(s)} \zeta(\bs)\,, \quad q\to 1 \,, 
$$
which defines a homomorphism $\grqMZV\to \MZV$, where 
$\grqMZV$ is the associated graded algebra for the weight 
filtration of $\qMZV$ and $\# \textup{odd}(\bs)$ is the 
number of odd terms in $\bs$. This homomorphism has a large kernel, 
already for $\QM \to \Q[\pi^2]$. Since $\zeta(\bs)>0$, we 
see that $Z(\bs)$ is not contained in any smaller weight 
filtration subspace. 

\subsection{}

We propose the following 

\begin{Conjecture}
The algebra $\qMZV$ is spanned by $Z(\bs)$ with $2 \le s_i \le 5$, 
$\Z/2$-graded by weight, and stable under
the operator $q \frac{d}{dq}$ that increases the weight by $2$\,.
The Hilbert series of the graded algebra $\grqMZV$ equals 
\begin{equation}
\sum_k t^k \dim_{\Q} \grqMZV_k =
\frac{1}{1-t^2-t^3-t^4-t^5+t^8+t^9+t^{10}+t^{11}+t^{12}} \,.
\label{Pser}
\end{equation}
\end{Conjecture}

The first statement here may be compared to 
a conjecture of Hoffman \cite{Hoff}, proven by F.~Brown \cite{Brown}, that says
$\zeta(\bs)$ with $s_i \in \{2,3\}$ span $\MZV$. However, 
while such $\zeta(\bs)$ are conjectured to be a basis of $\MZV$, 
the Hilbert series \eqref{Pser} predicts relations among 
$Z(\bs)$ with $2 \le s_i \le 5$\,. 

\section{Hilbert schemes of surfaces}\label{s2}

\subsection{}

Let $S$ be a nonsingular quasi-projective surface. The Hilbert 
scheme $\Hilb(S,n)$ parametrizes $0$-dimensional 
subschemes $\fS\subset S$ of length $n$, see e.g.\ \cite{Goe,Lehn,Nak}
for an introduction. 
It is an irreducible 
nonsingular quasi-projective variety of dimension $2n$. 

The geometry of $\Hilb(S,n)$, and in particular, the characteristic 
numbers of natural vector bundles on it are of 
great interest to algebraic geometers and mathematical physicists. 
 These characteristic numbers may be defined if 
$S$ is proper, or if there is a torus action on $S$ with proper 
fixed locus. In the latter case,  they take values in localized 
$G$-equivariant cohomology of a point. Here we fix
a connected reductive algebraic group $G$ that acts on $S$
so that the fixed-point set of its maximal torus is proper. 
It is convenient not to assume this action faithful. 

\subsection{}

A line bundle $\cL$ on $S$ defines a rank $n$ vector bundle on 
$\Hilb(S,n)$ with fiber 
$$
\cL^{[n]} \Big|_{\fS} = H^0(\cO_\fS \otimes \cL) \,. 
$$
These vector bundles are called tautological. 

The tangent bundle $\cT$ to the Hilbert scheme is described by 
$$
\cT\Big|_{\fS} = \chi(\cO_S) - \chi(\cI_\fS,\cI_\fS) \,, 
$$
where $\cI_\fS$ is the ideal sheaf of $\fS$ and 
$$
\chi(A,B) = \sum (-1)^i \Ext^i (A,B) \,. 
$$
{}From this, and the Grothendieck-Riemann-Roch theorem, it is 
possible to express the characteristic classes of $\cT$ in terms of 
those of $\cO^{[n]}$. For our purposes, however, 
such reduction appears quite impractical and it is 
perhaps best to keep the characteristic classes of $\cT$ separate from 
those of tautological bundles. 

Further, the tangent bundle $\cT$ may be twisted by a line bundle 
$\cM\in \Pic(S)$ as follows 
$$
\cT(\cM)\Big|_{\fS} = \chi(\cM) - \chi(\cI_\fS,\cI_\fS\otimes \cM)
\,. 
$$

In particular, if $\cM$ is a pure $G$-character  then $\cT(\cM) = 
\cT \otimes \cM$ and then $\cM$ would be called the mass 
of the adjoint matter in Nekrasov theory \cite{NN}. 

\subsection{}

Fix $\cL$, $\cM$, a characteristic class $f$ and form 
the following generating function 
\begin{equation}
\langle f \rangle = \sum_n q^n \int_{\Hilb(S,n)} 
f(\cL^{[n]}) \, \Euler(\cT(\cM))  \,.\label{<f>} 
\end{equation}
By construction 
$$
\langle f \rangle \in H^*_G(\pt)_\textup{loc}[[q]]\,. 
$$
Since $\cM$ may be always additionally twisted by a character, 
this is a generating function for the integrals of arbitrary 
characteristic class of $\cL^{[n]}$ against a Chern class of
$\cT(\cM)$. Precisely these combinations often come up 
in practice, see e.g.\ \cite{KST}, but not always packaged in this 
particular form. 

Note that for a nontrivial group $G$, the ring
$H^*_G(\pt)_\textup{loc}$ is nonzero in all degrees, so the
degrees of characteristic classes in \eqref{<f>} don't need to sum up to the 
dimension of the Hilbert scheme.

\subsection{} 

It was shown in \cite{CO} that 
$$
\langle 1 \rangle = \prod_{n>0} (1-q^n)^{\bd}\,, 
\quad 
\bd= - \int_S c_2(TS \otimes \cM) \,,
$$
where $TS$ is the tangent bundle of $S$.  We define
$$
\langle f \rangle' = \langle f \rangle \big/ \langle 1 \rangle 
$$
and define the connected generating functions by 
\begin{align*}
  \label{eq:1}
  \lla \ch_a \ch_b \rra = &\, \langle \ch_a \ch_b \rangle' - 
\langle \ch_a \rangle' \langle \ch_b \rangle'\\
\lla \ch_a \ch_b \ch_c \rra = &\, \langle \ch_a \ch_b \ch_c \rangle' - 
\langle \ch_a \ch_b \rangle' \langle \ch_c \rangle'\\
& - 
\langle \ch_a \ch_c \rangle' \langle \ch_b \rangle'
- 
\langle \ch_b \ch_c \rangle' \langle \ch_a \rangle'
+ 2\langle \ch_a \rangle' \langle  \ch_b \rangle' \langle \ch_c \rangle'\,,
\end{align*}
et cetera. Here $\ch_a$ are the components of the Chern character. 

The connected generating functions satisfy better integrality: 
$$
\lla f \rra \in \Bbbk[[q]] \,, 
$$
where $\Bbbk$ is the image of the map 
$$
H^*_G(S) \owns \gamma \mapsto \int_S \gamma \cup 
c_2(TS \otimes \cM) 
\in H^*_G(\pt)_\textup{loc} \,. 
$$

\subsection{}

These $q$-series are the subject of the following 

\begin{Conjecture}\label{c2} 
The series $\langle f \rangle'$ is a multiple $q$-zeta value of the 
same weight as $f$, where 
$$
\weight \ch_k = k+2 \,. 
$$
In the case $\cL  = \cK^{1/2}_S$, it is also of the same parity as the 
weight of $f$. 
\end{Conjecture}

Note that 
$$
\langle f \cdot \ch_0 \rangle'  = q\frac{d}{dq} 
\big[\langle f \rangle' + \ln \langle 1 \rangle \big]\,,
$$
which is why we were interested in the action of the operator 
$q\frac{d}{dq}$ on $\qMZV$. 

\subsection{}
Among the conjectures presented in this paper, Conjecture \ref{c2} 
appears the most accessible, perhaps using the techniques 
developed in \cite{CO}. In fact, prompted by this conjecture, it 
was already shown in \cite{CO} that 
$$
\langle c_1(\cO^{[n]}) \rangle' = 
\tfrac12\left(Z(2)-Z(3)\right) \, \int_S (c_1 c_2 - c_3) (TS \oplus
\cM)\,. 
$$ 
The main result of \cite{CO} computes \eqref{<f>} as the trace over 
the Fock space of a product of a certain vertex operator and a 
certain integral of motion of the second quantized trigonometric 
Calogero system. Formulas for the latter in terms of bosonic 
operators may be derived systematically using, for example, 
the formulas of \cite{Smir} or, alternatively, using 
many other approaches 
developed in the literature.  The trace can then be explicitly 
computed as a multiple $q$-series, similar to the series
$$
\sum_{k,l>0} \frac{q^{k+l}}{(1-q^k)(1-q^l)(1-q^{k+l})}\,,
$$
which was studied in \cite{CO}. The problem is thus reduced 
to showing that certain rather concrete series lie in $\qMZV$. 
It is hard to know without trying how difficult this task would be.

\section{Hilbert schemes of threefolds}

\subsection{}

Now let $X$ be a nonsingular quasiprojective threefold and 
consider its Hilbert scheme of points $\Hilb(X,n)$. This is 
a quite singular reducible scheme, however, as 
one of the simplest moduli space in Donaldson-Thomas theory 
\cite{Thom}  it has a perfect obstruction theory (here, of 
virtual dimension $0$) and the corresponding $0$-dimensional 
virtual fundamental class. 

In parallel to \eqref{<f>} we define 
\begin{equation}
\langle f \rangle_\td = \sum_n (-q)^n \int_{[\Hilb(X,n)]_\vir}
f(\cL^{[n]}) \,,  \label{f03D} 
\end{equation}
where $\cL^{[n]}$ is a rank $n$ bundle on $\Hilb(X,n)$ defined
in the same way as before. 

For example, $X$ could be the total space of a line bundle $\cM$ 
over a surface $S$ and then $\Hilb(S,n)$ is one of the component 
of $\Hilb(X,n)^\Ct$, where $\Ct$ acts by scaling $\cM$. The 
contribution of this component to \eqref{f3D} is closely
related to the integrals considered in Section \ref{s2}. However, 
the precise packaging of those integrals in the generating 
function is different, as will be the functional nature of the 
series \eqref{f03D}. 

\subsection{}

The evaluation 
\begin{equation}
\langle 1 \rangle_\td  = \prod_{n>0}(1-q^n)^{n \, \bdt} \,, 
\quad \bdt = \int_X (c_1 c_2 - c_3)\,, \label{13D}
\end{equation}
was checked for $X=\C^3$ in \cite{mnop} and  conjectured in general. 
It was proven for all $X$ in \cite{LP}, another proof was announced
by Jun Li.  As before, we set 
\begin{equation}
\langle f \rangle'_\td = \langle f \rangle_\td \big/ \langle 1
\rangle_\td \label{f3D}
\end{equation}
and define connected generating functions as in Section \ref{s2}. 
We also define a sequence of $H^*_G(\pt)$-submodules 
$$
H^*_G(\pt)_\textup{loc} \supset \Bbbk_1 \supset \Bbbk_2 \supset \dots
$$
as the images of the maps 
$$
H^*_G(X) \owns \gamma \mapsto \int_X \gamma \cup 
(c_1 c_2 -c_3)^d\,, \quad q=1,2,\dots 
$$

\subsection{}

Conjecturally, the generating functions \eqref{f3D} belong to the 
algebra 
\begin{equation}
\oqZ = \Q\left[\left(q\frac{d}{dq}\right)^{\! l} Z(2k+1) \right]_{k\ge 1, l
  \ge 0} \subset \qMZV
\label{oqZ}
\end{equation}
generated by the odd $q$-zeta values and their derivatives. This is 
the $q$-analog of the subalgebra 
$$
\oZ = \Q[\zeta(3),\zeta(5),\zeta(7), \dots] \subset \MZV\,,
$$
which, by standard transcendence conjectures, is a free commutative
$\Q$-algebra on its generators. A parallel conjecture for \eqref{oqZ} 
would be that it is also a free commutative algebra on its
generators. This allows us to define a grading 
by 3D weight and depth 
by setting 
\begin{align}
\weight_\td \left(q\frac{d}{dq}\right)^{\! l} Z(2k+1) &= 2k+1+l
\,,\label{wt3d}\\
\depth \left(q\frac{d}{dq}\right)^{\! l} Z(2k+1) &= 1 \notag
  \end{align}
on generators. 

\subsection{}

\begin{Conjecture}
Connected generating functions satisfy
\begin{equation}
  \label{eq:4}
 \lla \prod \ch_{k_i} \rra \in 
\sum_d \Bbbk_d \otimes_{\Q} \oqZ_{\le 2 + \sum (k_i+1),\ge d}\,,
\end{equation}
where first subscript denotes an upper bound on the 3D weight 
\eqref{wt3d} while the second denotes the depth $d$. 
\end{Conjecture}

\noindent
Presumably, the techniques used to prove \eqref{13D} will 
similarly reduce the general case of this conjecture to the 
case $X=\C^3$. 

\subsection{}

For $X=\C^3$, any given instance of this conjecture may be 
attacked, at least in principle, by expanding in a series in
$$
t_{12} = t_1 + t_2 \in H^2_{GL(3)}(\pt)\,, 
$$
where $t_i$ are the weights of the coordinates of $\C^3$. 
For connected functions, we have
$$
 \lla \prod \ch_{k_i} \rra \in \bdt \cdot H^{\sum k_i}_{GL(3)}(\pt)[[q]]\,, 
\quad \bdt = \tfrac{(t_1+t_2)(t_1+t_3)(t_2+t_3)}{t_1 t_2 t_3}\,, 
$$
and so the expansion in $t_{12}$ is a finite expansion. 

Let $\pi$ is a 3-dimensional partition, presented as piles of 
boxes $\bx=(i,j,k)$ of height $\pi_{ij}$, that is
$$
1 \le k \le \pi_{ij}\,, 
$$
placed over the squares
$\square = (i,j)$ of a 2-dimensional partition $\lambda$ lying in the
$(x_1,x_2)$-plane. The level sets of the function $\pi$ define a 
a decomposition of $\lambda$ into skew diagrams 
$$
\varnothing \subset \dots \subset \lambda'' \subset \lambda' 
\subset \lambda \,. 
$$
For a skew diagram, we define its rank as the minimal number 
of rim hooks needed to decompose it, see Section 4.5 of
\cite{OP}. We define the 
total rank of $\pi$ by 
$$
\rho(\pi) = \sum 
\textup{rank}(\lambda^{(i)}/\lambda^{(i+1)})
$$
It can then 
be shown \cite{OP} that 
$$
\textup{order}_{t_{12}}  \, \textup{contribution}(\pi) = \rho(\pi) 
$$
where order is the order of vanishing along $t_{12}=0$ and the 
contribution of $\pi$ is the weight at $\pi$ of the virtual 
fundamental cycle in the equivariant localization formula. 

Terms of small degree in $t_{12}$ thus come from 3-dimensional 
partitions of small total rank, and it may be possible to analyze
them directly. 

\subsection{}

For example, consider the computation of 
$$
\langle \ch_0 \rangle'  \in \bdt \cdot \Q[[q]] \,, 
$$
for which it is enough to compute the linear term in $t_{12}$. 
A 3-dimensional partition 
$\pi$ has total rank $1$ if and only if 
\begin{enumerate}
\item partition $\lambda$ is a hook, and 
\item the function $\pi_{ij}$ is a constant, say, $c$. 
\end{enumerate}
Denoting $n=|\pi|$ the size of $n$, we see that $c$ is a divisor 
of $n$ and, for given $c$, there are exactly $n/c$ choices of 
$\lambda$. One further checks that 
$$
\left. \left( \frac{t_3}{t_{12}} \, \textup{contribution}(\pi)\right)
\right|_{t_{12}=0,t_3=0} = \frac{ (-1)^{n+1}}{c} \,.
$$
Since $\ch_0 = n$, we obtain 
$$
\langle \ch_0 \rangle' = - \bdt \sum_n q^n \sum_{c|n} (n/c)^2  
= - \bdt Z(3) \,, 
$$
recovering \eqref{13D}. 

\subsection{}
It is interesting to find out to what extent the
Gromov-Witten/Donaldson-Thomas correspondence  of \cite{mnop} holds for 
the functions \eqref{f3D}. 

Object parallel to the connected 
functions $\lla \prod \ch_{k_i} \rra$ on the Gromov-Witten side are the following 
Hodge integrals 
\begin{equation}
\left\langle \prod_{i=1}^n \tau_{k_i+1} \right\rangle = 
\sum_{g\ge 0} u^{2g-2} \int_{\Mgn} 
\Lambda(t_1) \Lambda(t_2) \Lambda(t_3) 
\prod_{i=1}^n \psi_i^{k_i+1} \label{Hodge}
\end{equation}
where $\Mgn$ is the Deligne-Mumford moduli space of 
stable $n$-pointed genus $g$ curves $C$, 
$$
\Lambda(t) = t^g-t^{g-1} \lambda_1 + \dots\pm \lambda_g \,, 
\quad \lambda_i = c_i (H^0(\omega_C))\,, 
$$
is, up to normalization, the Chern polynomial of the Hodge bundle, and 
$\psi_i \in H^2(\Mgn)$ is 1st Chern class of the line bundle formed 
by the cotangent line at the $i$th marked
point. 

\subsection{}

A correspondence would involve knowing the top and the 
bottom arrow in the following diagram 
$$
\xymatrix{
\Q[\ch_k]  \ar@{->}[rr]^\Phi \ar@{->}[d]_{\lla \,\,\, \rra} && \Q[i][\tau_k] 
 \ar@{->}[d]^{\langle \,\,\, \rangle} \\
 \Bbbk \otimes \oqZ   \ar@{->}[rr]^{\phi} 
&& \Bbbk \otimes \Q[i][u^{-1},u]] 
}
$$
in which $\phi$, but not $\Phi$, should be an algebra 
homomorphism. For the top map $\Phi$ there are certain 
general principles stated in \cite{mnop} and a also concrete explicit 
proposal \cite{OOP}. In the usual GW/DT story, the bottom map is 
given by the expansion of a rational function of $q$ in a Laurent 
series in $u$ via the substitution
 $q=e^{iu}$. This cannot literally work in the present
case for at least two reasons. 

First, certain factors of $\frac{1}{2}$ appear consistently in the 
treatment of the GW/DT correspondence for Hilbert scheme of points 
(as opposed to curves of nonzero degree). This concerns even the 
series \eqref{13D}. Second, the asymptotic expansion of the odd 
zeta values contains transcendental terms, namely 
\begin{equation*}
  \label{EM}
  Z(2k+1)
\sim 
\frac{(2k)!\, \zeta(2k+1)}{(-\log q)^{2k+1}} - \sum_{n=0}^\infty
\frac{B_{2n+2k} \, B_{2n}}{ (2n)!\, (2n+2k)} \, (-\log q)^{2n-1} \,,
\end{equation*}
as $q\uparrow 1$. If one sets
$$
\phi\left( Z(2k+1)\right) = - \boldsymbol{\frac12} \sum_{n=0}^\infty
\frac{B_{2n+2k} \, B_{2n}}{ (2n)!\, (2n+2k)} \, (-iu )^{2n-1}
$$
then for small examples one can find an agreement between 
the Donaldson-Thomas and Gromov-Witten computations. However, 
I am not convinced that it would work in general. 

Instead, I would like to pose the determination of the functional 
nature of the series \eqref{Hodge} as an open problem, the solution 
to which may very well require a generalization of the universe of 
asymptotic expansions of $q$-zeta values. 

\subsection{Acknowledgments} 

I am very grateful to Erik Carlsson, Alexander Goncharov, Davesh Maulik, 
Nikita Nekrasov, Alexei Oblomkov, and Rahul Pandharipande for 
discussions. I thank NSF for its support in the framework of FRG 1159416.

\newpage 

\noindent 
Andrei Okounkov\\
Department of Mathematics, Columbia University\\
New York, NY 10027, U.S.A.\\

\vspace{-12 pt}

\noindent 
Institute for Problems of Information Transmission\\
Bolshoy Karetny 19, Moscow 127994, Russia\\

\vspace{-12 pt}

\noindent 
Laboratory of Representation
Theory and Mathematical Physics \\
Higher School of Economics \\ 
Myasnitskaya 20, Moscow 101000, Russia

\end{document}